\documentclass[dblspace,leqno]{aoscm}
\usepackage{graphics}
\usepackage{natbib}
\bibpunct{(}{)}{;}{x}{,}{,}
\usepackage[colorlinks=true,citecolor=blue]{hyperref}
\usepackage{amsmath}
\usepackage{amssymb}
\usepackage{fancybox}
\usepackage{fancyheadings}

\def\theequation{\arabic{equation}} 
\def\thefigure{\arabic{figure}} 

\def\Cov{\mathop{\rm cov}\nolimits}

\newcommand{\ip}[2]{\left <{#1},{#2}\right >}

\renewcommand{\to}{\rightarrow}

\newcommand{\mbP}{\mathbb{P}}
\newcommand{\mbV}{\mathbb{V}}

\newcounter{assumption}

\newtheorem{theorem}{Theorem}
\newtheorem{lemma}{Lemma}

\newtheorem{proposition}{Proposition}

\newcounter{remark}

\newcounter{example}

\newcommand{\bc}{\begin{center}}
\newcommand{\ec}{\end{center}}
\newcommand{\be}{\begin{equation}}
\newcommand{\ee}{\end{equation}}
\newcommand{\bea}{\begin{eqnarray}}
\newcommand{\eea}{\end{eqnarray}}
\newcommand{\bes}{\begin{eqnarray*}}
\newcommand{\ees}{\end{eqnarray*}}
\newcommand{\bed}{\begin{description}}
\newcommand{\eed}{\end{description}}
\newcommand{\bee}{\begin{enumerate}}
\newcommand{\eee}{\end{enumerate}}
\newcommand{{\mbt}}{{\mbox{\boldmath $\beta$}} }
\newcommand{{\wmbt}}{\widehat{\mbt}}
\newcommand{{\mbg}}{{\mbox{\boldmath $\gamma$}} }
\newcommand{{\wtmbg}}{\widetilde{\mbg}}
\newcommand{{\wmbg}}{\widehat{\mbg}}
\newcommand{{\mbr}}{{\mbox{\boldmath $\theta$}} }
\newcommand{{\wmbr}}{\widehat{\mbr}}
\newcommand{{\wtmbr}}{\widetilde{\mbr}}
\newcommand{{\mbx}}{{\mbox{\boldmath $\xi$}} }
\newcommand{{\mbe}}{{\mbox{\boldmath $\eta$}} }
\newcommand{{\mve}}{{\mbox{\boldmath $\epsilon$}} }
\newcommand{{\mvv}}{{\mbox{\boldmath $\vartheta$}} }
\newcommand{\wmbx}{\widehat{\mbx}}

\newcommand{{\mba}}{{\mbox{\boldmath $\alpha$}} }

\def\b#1{\mbox{\boldmath$#1$}}

\newcommand{{\mbn}}{ {\bf u} }
\newcommand{{\mbm}}{{\mbox{\boldmath $\mu$}} }
\newcommand{{\mbl}}{{\mbox{\boldmath $\lambda$}} }
\newcommand{{\mbv}}{{\mbox{\boldmath $\varpi$}} }
\newcommand{{\mbrh}}{{\mbox{\boldmath $\rho$}} }

\newcommand{{\Om}}{{\mbox{\boldmath $\Omega$}} }

\newcommand{\mV}{\mathcal{V}}
\newcommand{\mH}{\mathcal{H}}
\newcommand{\mK}{\mathcal{K}}
\newcommand{\mM}{\mathcal{M}}
\newcommand{\mB}{\mathcal{B}}
\newcommand{\mT}{\mathcal{T}}
\newcommand{\fT}{\mathfrak{T}}
\newcommand{\mN}{\mathcal{N}}
\newcommand{\mG}{\mathcal{G}}
\newcommand{\mS}{\mathcal{S}}
\newcommand{\mC}{\mathfrak{C}}
\newcommand{\mP}{\mathcal{P}}

\newcommand{\mQ}{\mathcal{Q}}

\newcommand{\mE}{\mathbb{E}}
\newcommand{\mbT}{{\bf T}}
\newcommand{\mbQ}{{\bf S}}

\newcommand{\bbL}{\mathcal{L}}

\newcommand{\bB}{{\bf B}}
\newcommand{\bQ}{{\bf Q}}

\newcommand{\g}{{\rm g}}
\newcommand{\bg}{{\b \g} }

\renewcommand{\a}{{\rm a}}
\newcommand{\bb}{{\rm b}}

\newcommand{\bA}{{\bf A}}
\newcommand{\bR}{{\bf R}}
\newcommand{\bD}{{\bf D}}
\newcommand{\bL}{{\bf L}}
\newcommand{\bW}{{\bf W}}
\newcommand{\bY}{{\bf Y}}
\newcommand{\bZ}{{\bf Z}}
\newcommand{\wZ}{\widehat{\bZ}}

\newcommand{\bU}{{\bf U}}

\newcommand{\bX}{{\bf X}}
\renewcommand{\S}{{\rm S}}

\newcommand{\bh}{{\bf h}}
\newcommand{\bC}{{\bf C}}
\newcommand{\wC}{\widehat{\bC}}
\newcommand{\bP}{{\bf P}}
\newcommand{\bM}{{\bf M}}
\newcommand{\bI}{{\bf I}}
\newcommand{\bJ}{{\bf J}}
\newcommand{\bH}{{\bf H}}

\newcommand{\mbF}{{\mathbf F}}
\newcommand{\bS}{{\mathbf \Sigma}}

\NONUMBIB
\begin{document}

\SPECFNSYMBOL{1}{2}{3}{}{}{}{}{}{}%

\AOSMAKETITLE

%\AOSyr{2006}
%\AOSvol{00}
%\AOSno{00}
%\AOSpp{000--000}
%\AOSReceived{Received ***}
\AOSAMS{Primary 62F30; Secondary 62F03, 62F05.}
\AOSKeywords{Convex analysis; Gauss-Bonnet theorem;
Generalized quasi-score test; Likelihood ratio test; Local
alternatives; Mixed volumes; Monotone regression; Order restricted
hypothesis; Volume-of-tube formula; Weighted chi-squared.}

\AOStitle{Inference Under Convex Cone Alternatives for Correlated
Data} 
\AOSauthor{Ramani S. Pilla\footnote{Research supported in part by the
National Science Foundation grant DMS 02-39053 and the Office of Naval
Research grants N00014-02-1-0316 and N00014-04-1-0481.} }

\AOSaffil{Case Western Reserve University}
\AOSlrh{RAMANI S. PILLA}
\AOSrrh{INFERENCE UNDER CONVEX CONE ALTERNATIVES}
\AOSAbstract{ 
In this research, inferential theory for hypothesis testing under
general convex cone alternatives for correlated data is
developed. While there exists extensive theory for hypothesis testing
under smooth cone alternatives with independent observations,
extension to correlated data under general convex cone alternatives
remains an open problem. This long-pending problem is addressed by (1)
establishing that a {\em generalized quasi-score} statistic is
asymptotically equivalent to the squared length of the projection of
the standard Gaussian vector onto the convex cone and (2) showing that
the asymptotic null distribution of the test statistic is a weighted
chi-squared distribution, where the weights are {\em mixed volumes} of
the convex cone and its polar cone. Explicit expressions for these
weights are derived using the volume-of-tube formula around a convex
manifold in the unit sphere. Furthermore, an asymptotic lower bound is
constructed for the power of the generalized quasi-score test under a
sequence of local alternatives in the convex cone. Applications to
testing under order restricted alternatives for correlated data are
illustrated. }
\maketitle
\section{Introduction}
\label{sec-intro}
Correlated or longitudinal data arise in many areas of science when a
response is measured at repeated instances on a set of subjects. It is
assumed that the measurements on different subjects are independent,
while those on individual subjects are correlated with an unknown
correlation structure \citep{diggle:94}. In this research, inferential
theory is developed for hypothesis testing under general convex cone
alternatives for correlated data using the formula for the volume of a
tube around a manifold (curve, surface, etc.) on the surface of the
unit sphere in an $r$-dimensional Euclidean space $\Re^r$
\citep{hotelling:39,weyl:39,naiman:90}. Testing for order
restricted parameters in correlated data and testing for a monotone
regression become special cases of this general problem. Often,
interest lies in detecting an order among treatment effects, while
simultaneously modeling relationships with regression
parameters. There exists extensive theory for hypothesis testing under
ordered alternatives with independent observations
\citep{barlow:72,robert:88,silva:04}, including smooth cone alternatives
\citep{take:97}. However, extension of the theory to correlated data
remains an open problem.

\subsection{\em Formulation of the testing problem} 
Let $Y_{ij}$ be the response measured at the $j$th ($j = 1, \ldots,
n_{i}$) time point on the $i$th ($i = 1, \ldots, N$) subject. Let
$\bY_{\!i} = (Y_{i1}, \ldots, Y_{i n_{i}})^{\, T}$ be an
$n_i$-dimensional vector of response variables.  The mean of $Y_{ij}$
is related to the $r$-dimensional vector of covariates $\bX_{ij}$
corresponding to the $r$-dimensional parameter vector $\mbg$ via a
generalized linear model 
\bea
  \label{eq:link}
  \mE(Y_{ij}) := h \big(\bX_{ij}^{\, T} \, \mbg \big),
\eea
where $h(\cdot)$ is the inverse of a link function. We assume that the
true distribution is unique and all expectations are taken with
respect to the true probability measure $P$.

The goal is to test the general hypothesis
\bea
   \label{eq:hyp0}
   \mH_0\!: \mbg \in \mV  \quad \mbox{against} \quad \mH_1\!: \mbg
   \in \mV \oplus \mC, \mbg \notin \mV,
\eea
where $\mV$ is an arbitrary finite dimensional vector space of
$\Re^r$, $r := \mbox{dim}(\mV \oplus \mC)$, $\mC$ is a closed convex
cone with a non-empty interior in $\Re^r$ and $\oplus$ denotes the
{\em direct} or {\em Kronecker sum}.  Without loss of generality, it
is assumed that $\mC \subset \mV^{\perp}$, the orthogonal complement
of $\mV$. Under $\mH_1$, $\mbox{dim}(\mbg) = r$; whereas under
$\mH_0$, $\mbox{dim}(\mV) < r$ due to certain constraints imposed on
the parameters in $\mV$. 

Seminal work of \cite{take:97} has established a solution for the
problem of testing a simple null hypothesis regarding the multivariate
Gaussian mean vector $\mbl$ against an arbitrary convex cone
alternative for independent observations. In particular, they derived
the asymptotic null distribution of the likelihood ratio test
statistic (LRT) for testing
\bea
  \label{eq:hyp1}
  \mH_0^{\, 1}\!: \mbl = {\bf 0} \quad \mbox{against} \quad \mH_1^{\,
  1}\!: \mbl \in \mK, 
\eea
where $\mK$ is a closed convex cone of dimension $d$ with a nonempty
interior in $\Re^r$ ($r \geq d$), using the techniques of convex
analysis. 

\subsection{\em Main results and organization of the article} 
The goals of this research include the following.
\bee
\vspace{-0.15in}
\item In Section \ref{sec-qif}, we derive large-sample properties
of the quadratic inference functions, extensions of the generalized
method of moments \citep{hansen:82}, that are required for the
development of inferential theory with correlated data.
\vspace{-0.15in}
\item We derive a ``generalized quasi-score'' (GQS) statistic for
the testing problem (\ref{eq:hyp0}) for correlated data in Section
\ref{sec-gcone}. Furthermore, it is established that the asymptotic
null distribution of the GQS statistic, to appropriate statistical
order, is equivalent to finding the limiting distribution of the
squared length of projection of a standard Gaussian vector onto the
convex cone $\mK$ (see Theorem \ref{thm-proj}).
\vspace{-0.15in}
\item In Section \ref{sec-gcone}, we also establish that the
asymptotic null distribution of the GQS statistic is a weighted
chi-squared distribution, where the weights are {\em mixed volumes} of
$\mK$ \citep{take:97} and its polar cone (see Theorem
\ref{them-adist}). Deriving computable expressions for the weights in
the asymptotic null distribution of the test statistic is a tedious
and difficult process even for independent data \citep{take:97}. Only
for special cases, weights are known explicitly or can be determined
numerically. 
\vspace{-0.15in}
\item In Section \ref{sec-wts}, we express the asymptotic null
distribution of the GQS statistic in terms of certain geometric
constants of the volume-of-tube formula 
\citep{hotelling:39,weyl:39,adler:81,naiman:90} for Gaussian random
fields \citep{siegmund:95, worsley2:95,worsley1:95,worsley:96} around
a {\em convex manifold} on the surface of the unit sphere (see Theorem
\ref{th:sig3d}). We derive explicit expressions (in suitable forms for
computation) for these geometric constants by representing them as
integrals over appropriate parts of the manifold. 
\vspace{-0.15in}
\item We derive an asymptotic lower bound for the power of the GQS
test under a sequence of local alternatives in $\mK$ in Section
\ref{sec-lalt} (see Theorem \ref{thm-lalt}). This lower bound
demonstrates that the test under restricted alternatives is more
powerful than the corresponding one under unrestricted
alternatives. To the best of the author's knowledge, no such lower
bound has been derived in the literature even for independent data.
\vspace{-0.15in}
\eee
The article concludes with a discussion in Section \ref{sec-discuss}. 

\section{Large-Sample Properties of the Inference Functions}
\label{sec-qif}
In this section, we present the large-sample properties of the
inference functions that are required for later theoretical
development of our general testing problem. See \cite{pilla:05} for
technical details and derivations of the results presented in this 
section. 

\cite{hansen:82} proposed the {\em generalized method of moments}
(GMMs) for estimating the vector of regression parameters $\mbt \in
\mB$ from a set of score functions, where the dimension of the score
function exceeds that of $\mbt$. He established that, under certain
regularity conditions, the GMM estimator is consistent, asymptotically
Gaussian, and asymptotically efficient. \cite{qu:00} extended the GMMs
to create a clever approach called the ``quadratic inference
function'' (QIF) that implicitly estimates the underlying correlation
structure for the analysis of longitudinal data. Their main idea was
to assume that the inverse of the working correlation matrix, denoted
by $\bR^{-1}(\mba)$, is a linear combination of several pre-specified
basis matrices. That is, $\bR^{-1}(\mba) = \sum_{l = 1}^{s}
\alpha_l \, \bM_l$, where $\alpha_1, \ldots, \alpha_{s}$ are
unknown constants, $\bM_1$ is the identity matrix of an appropriate
dimension and $\bM_l$ $(l = 2, \ldots, {s})$ are pre-specified
symmetric matrices with elements taking either 0 or 1 for the commonly
employed working correlation structures such as exchangeable, AR-1
etc.

\subsection{\em Properties of extended score functions}
\label{sec-CN}
For mathematical exposition, we assume that each subject is observed
at a common set of times $j = 1, \ldots, n$. Let $\bh_{i} = \left[ 
h \big(\bX_{i1}^{\, T} \, \mbg \big), \ldots, h \big(\bX_{in}^{\, T}
\, \mbg \big) \right]^{\, T}$, where $h \big(\bX_{ij}^{\, T} \, \mbg
\big)$ is the inverse of a link function and the operator $\nabla$
denotes partial derivative with respect to the elements of $\mbg$;
therefore, $\nabla \, \bh_{i}$ is the $(n \times r)$ matrix
$\left(\partial \, \bh_{i}/\partial \, \gamma_1, \ldots, \partial \,
\bh_{i}/\partial \, \gamma_r \right)$ for each $i = 1, \ldots, N$. 

The coefficients $\alpha_1, \ldots, \alpha_{s}$ in $\bR^{-1}(\mba)$
are treated as nuisance parameters to create the set of
subject-specific {\em basic score functions} as 
\bea
  \label{eq:escore1}
  \bg_{i}(\mbg) &:=& \left[ \begin{array}{c} \nabla \bh_{i}^{T}
  \; \bA_{i}^{-1/2} \; \bM_1 \; \bA_{i}^{-1/2} \;
  \left( \bY_{i} - \bh_{i} \right) \\ \vdots \\ \nabla \bh_{i}^{T}
  \; \bA_{i}^{-1/2} \; \bM_{s} \; \bA_{i}^{-1/2}
  \, \left( \bY_{i} - \bh_{i} \right) \end{array} \right] \quad
  \mbox{for} \quad i = 1, \ldots, N,
\eea
where $\bA_{i}$ is the diagonal matrix of marginal covariance of
$\bY_i$ for the $i$th subject.

Define the vector of {\em extended score functions} for all
subjects as $\overline{\bg}_N(\mbg) := N^{-1} \, \sum_{i = 1}^N
\bg_i(\mbg)$. Note that $\mbox{dim}[\overline{\bg}_{N}(\mbg)] =
r \, s > r = \mbox{dim}(\mbg)$. The extended score vector
$\overline{\bg}_{N}(\mbg)$ satisfies the mean zero assumption
$\mE_{\mbg}[\overline{\bg}_N(\mbg)] = 0$, where the expectation
operator is with respect to the true but unknown distribution of the
response matrix $\bY$. These estimating equations can be combined 
optimally using the GMM \citep{hansen:82}. 

Let $\bS_{\mbg_{0}}(\mbg)$ be the true covariance matrix of
$\bg_{1}(\mbg)$, an $s$-dimensional vector of extended score functions
defined in (\ref{eq:escore1}). We require the following design
assumptions for deriving the asymptotic theory. \\
\begin{assumption}
\label{as0}
  The pairs $(\bY_{i}, \bX_{i}^T)$ for $i = 1, \ldots, N$, where
  $\bX_{i} = (\bX_{i1}, \ldots, \bX_{in})$ are $(r \times
  n)$-dimensional matrices, are an independent sample from an $[n
  \times (r + 1)]$-dimensional distribution $\mbF$.
\end{assumption} 
\\
\begin{assumption}
\label{as1}
The number of measurements $n_{i}$ on the $i$th subject is fixed
at $n_i = n$ for all $i = 1, \ldots, N$. 
\end{assumption}
\\
\begin{assumption}
\label{as2}
The $(rs \times rs)$-dimensional covariance matrix
$\bS_{\mbg_{0}}(\mbg_{0}) := \mE_{\mbg_{0}} [\bg_{1}(\mbg_0) \,
\bg_{1}^T(\mbg_0)]$ is strictly positive definite.
\end{assumption} \\
The independence part of the assumption \ref{as0} is between different
subjects (i.e., with respect to the index \(i\)). The elements of $\bX_i$
need not be independent of each other; hence, this assumption
incorporates both time-dependent as well as time-independent
covariates. Moreover, there exists a dependence of $\bY_i$ on $\bX_i$
through the link function given in (\ref{eq:link}). 

All throughout this article, $\mE_{\mbg_{0}}$ denotes an expectation
operator with respect to the true parameter vector $\mbg_{0}$ and all
expectations are assumed to be finite. Let $\wC_{N}(\mbg)$ be the
estimator of the {\em second moment matrix} of $\bg_1(\mbg)$ so that
$\wC_N(\mbg) := N^{-1} \sum_{i = 1}^N \bg_i(\mbg) \;
\bg_i^{T}(\mbg)$. If the mean zero assumption for
$\overline{\bg}_N(\mbg)$ holds, then $N^{-1} \, \wC_N(\mbg)$ estimates
the covariance of $\overline{\bg}_N(\mbg)$. 

\subsection{\em Fundamental results for the quadratic inference functions}
\label{sec-asymp}
The {\em quadratic inference function} (QIF) is defined as
\bea
  \label{eq:qif} 
  \mQ_N(\mbg) := N \, \overline{\bg}_N^T(\mbg) \;
  \wC_{N}^{-1}(\mbg) \; \overline{\bg}_N(\mbg).
\eea
If rank of $\wC_N(\mbg)$ is less than $rs$ or is singular, then the
inverse does not exist. However, any vector in the null space of
$\wC_N(\mbg)$ must be orthogonal to each of the subject-specific score
functions $\bg_i(\mbg) \; (i = 1, \ldots, N)$ and consequently to
$\overline{\bg}_N(\mbg)$. Therefore, one can replace
$\wC_N^{-1}(\mbg)$ by any generalized inverse such as the
Moore-Penrose generalized inverse. Our covariance estimator
$\wC_N(\mbg)$ differs from that of \cite{qu:00}, who define a
covariance $C_N$ with a factor of $N^{-2}$, and correspondingly omit
the factor of $N$ from the QIF in (\ref{eq:qif}). There are some
technical flaws in their results and hence we proceed carefully
without relying on those asymptotic results. \cite{pilla:05}
established precise large-sample results for the QIF along with
rigorous proofs. The QIF in (\ref{eq:qif}) is minimized under
unrestricted and restricted spaces to yield the estimators $\wmbg =
\underset{\mbg \, \in \, \Re^r}{\arg \, \min} \; \mQ_N(\mbg),
\widetilde{\mbg} = \underset{\mbg \, \in \, \mV \oplus \mC}{\arg \,
\min} \; \mQ_N(\mbg)$ and $\overline{\mbg} = \underset{\mbg \, \in \,
\mV}{\arg \, \min} \; \mQ_N(\mbg)$, respectively. These estimators can
be found using the {\em iterative reweighted generalized least
squares} (IRGLS) algorithm developed by \cite{loader:05}. The IRGLS
algorithm avoids the complexity of computing the second derivative
matrix of $\mQ_N(\mbg)$ required for employing the Newton-Raphson
algorithm recommended by \cite{qu:00}.

The proof of the next lemma essentially follows from p.\ 26 of
\cite{lee:96} and hence is omitted. 
\begin{lemma}
\label{thm11}
Under assumptions \ref{as0}--\ref{as2}, the QIF estimator $\wmbg
= \underset{\mbg \, \in \, \Re^r}{\arg \, \min} \; \mQ_N(\mbg)$ exists
uniquely and is strongly consistent. That is, $\wmbg
\overset{p}{\longrightarrow} \mbg_{0}$ as $N \to \infty$.  
\end{lemma}
We require the following regularity conditions for further theoretical
development. \\ 
\begin{assumption}
\label{as5} The parameter space of $\mbg$ denoted by $\mG \subset
\Re^r$ is compact. 
\end{assumption}
\\
\begin{assumption}
\label{as3} The parameter space of $\mbg$ is identifiable:
$\mE_{\mbg_{0}}[\bg_{1}(\mbg)] \neq 0$ if $\mbg \neq
\mbg_{0}$.  
\end{assumption} 
\\
\begin{assumption}
\label{as4} The true covariance matrix $\bS_{\mbg_{0}}(\mbg)$ is a
continuous function of $\mbg$.    
\end{assumption}
\\
\begin{assumption}
\label{as6}
The expectation $\mE_{\mbg_{0}}[\overline{\bg}_N(\mbg)]$ exists,
finite for all $\mbg \in \mG$ and continuous in $\mbg$.  
\end{assumption}
\\
\begin{assumption}
\label{as7}
The subject-specific score functions $\bg_{i}(\mbg)$ $(i =
1, \ldots, N)$ have uniformly continuous second-order partial
derivatives with respect to the elements of $\mbg$.
\end{assumption}
\\
\begin{assumption}
\label{as8}
The first-order partial derivatives of $\overline{\bg}_N(\mbg)$ and
$\wC_N(\mbg)$ have finite means and variances.
\end{assumption}

The importance of assumption \ref{as5} is that it enables us to invoke
Theorem 1 of \cite{rubin:56}; therefore, convergence statements in
this article are uniform for $\mbg$ in bounded sets.
\begin{lemma} 
\label{lemm1}
Under assumptions \ref{as0}--\ref{as4}, $\wC_N(\wmbg)
\overset{p}{\longrightarrow} \bS_{\mbg_{0}}(\mbg_0)$ as $N \to
\infty$.
\end{lemma} 
{\sc Proof.} Under assumptions \ref{as0}--\ref{as2} and 
by the strong law of large numbers, $\wC_N(\mbg)$ converges to its 
expected value, a non-degenerate limit, for a fixed $\mbg$. That is, 
\bea 
  \label{eq:cnhat}
  \wC_{N}(\mbg) \overset{p}{\longrightarrow} \; \mE[
  \bg_{1}(\mbg) \; \bg_{1}^T(\mbg)] \quad \mbox{as} \quad N
  \to \infty. 
\eea
Under the stated regularity conditions and Theorem 1 of
\cite{rubin:56}, uniform convergence holds if the compactness
assumption \ref{as5} holds. Hence, the claim (\ref{eq:cnhat}) holds
uniformly in $\mbg$. Lemma \ref{thm11} combined with the continuity of
the function $\wC_N(\mbg)$ ensures that $\wC_{N}(\wmbg)
\overset{p}{\longrightarrow} \; \mE[ \bg_{1}(\mbg_0) \;
\bg_{1}^T(\mbg_0)]$ as $N \to \infty$. The claim then follows from the
definition of $\bS_{\mbg_{0}}(\mbg_0)$. \hfill \rule{2mm}{2mm} 

Let $\bD(\mbg) := \mE_{\mbg}[ \nabla \bg_{1}(\mbg)]$, where $\nabla
\bg_{1}(\mbg) = \partial \, \bg_{1}(\mbg)/\partial \, \mbg$. From the
strong law of large numbers, $\nabla \overline{\bg}_N(\mbg)
\overset{p}{\longrightarrow} \bD(\mbg) \quad \mbox{as} \quad N \to
\infty$. This relation combined with Lemma \ref{lemm1}, enable us to
obtain the asymptotic covariance matrix of $\wmbg$. When there is no
ambiguity, we drop the subscript $\mbg_0$ and write
$\bS^{-1}(\mbg_{0})$ for the true covariance matrix of $\bg_1(\mbg)$
evaluated at $\mbg_0$. 

Let the estimated covariance of $\wmbg$ be defined as
\bea
  \label{covgam}
  \widehat{\Cov}(\wmbg) := \frac{1}{N} \left[ \nabla \,
  \overline{\bg}_N^T(\wmbg) \, \wC_N^{-1}(\wmbg) \,
  \nabla \overline{\bg}_N(\wmbg) \right]^{-1} \quad \mbox{as} \quad N
  \to \infty. 
\eea
\begin{lemma}
  \label{lemm2}
  Under assumptions \ref{as0}--\ref{as8}, $N \, \widehat{\Cov}(\wmbg)
  \overset{p}{\longrightarrow} \left[ \bD^{\, T}(\mbg_{0}) \;
  \bS^{-1}(\mbg_{0}) \; \bD(\mbg_{0}) \right]^{-1} =:
  \bJ^{-1}(\mbg_{0})$ as $N \to \infty$.  
\end{lemma}
The proof of the above lemma follows from the previous results. An
immediate consequence of Theorems 3.1 and 3.2 of \cite{hansen:82} is
the next result. The notation $\rightsquigarrow$ denotes convergence
in distribution. 
\begin{theorem}
  \label{thm12}
  Under assumptions \ref{as0}--\ref{as8}, $\sqrt{N} \; (\wmbg -
  \mbg_{0}) \rightsquigarrow N_{r} [{\bf 0}, \bJ^{-1}(\mbg_{0})]$ as
  $N \to \infty$.  
\end{theorem}
\begin{theorem}
\label{qnthm}
   Under assumptions \ref{as0}--\ref{as8}, $(2 \sqrt{N})^{-1} 
   \nabla \mQ_N(\mbg_0) \rightsquigarrow N_r [{\bf 0},
   \bJ(\mbg_{0})]$ as $N \to \infty$. 
\end{theorem}

\subsection{\em Testing under order restricted alternatives for correlated data}
\label{sec-order}
In the context of correlated data, comparing several treatments,
groups or populations with respect to their means, medians or location
parameters often arise in many areas of scientific applications. For
instance, one assumes that certain treatments are not worse than
another. 

The problem of testing under order restricted or constrained
hypothesis in longitudinal data becomes a special case of
(\ref{eq:hyp0}). Let $Y_{ijt}$ be the measurement taken at the $j$th
($j = 1, \ldots, n_{it}$) time point on the $i$th ($i = 1, \ldots,
n_t$) subject in the $t$th ($t = 1, \ldots, m$) treatment group. Let
$N = \sum_t n_t$. For mathematical exposition, we assume that $n_{it}
= n$ for all $(i, t)$ pairs. The mean of $Y_{ijt}$ is related to the
$p$-dimensional vector of covariates $\bX_{ijt}$, corresponding to the
$p$-dimensional parameter vector $\mbt_t$ for the $t$th group, via
$\mE(Y_{ijt}) := h
\left(\bX_{ijt}^{\, T} \, \mbt_t + \mu_{_t} \right)$, where $h(\cdot)$
is the inverse of a link function, $\mu_{_t}$ is the treatment effect
for the $t$th group. Hence, $\mbg = (\mbm^T, \mbt^T)^{\, T}$ with
$\mbm = (\mu_{_1}, \ldots, \mu_{_m})^{\, T}$ and $\mbt = (\mbt^{T}_1,
\ldots, \mbt^{T}_m)^{\, T}$. 

The order restricted hypothesis testing problem that is of interest
can be formulated as $\mH_0^{\, o}\!: \mbm \in \mV_0$ against
$\mH_1^{\, o}\!: \mbm \in \mC_0, \mbm \notin \mV_0$, where $\mV_0 =
\{\mbm\!: \mbm_1 = \cdots = \mbm_m\}$ and $\mC_0 = \{\mbm\!: \mbm_1
\geq \cdots \geq \mbm_m\}$ is a particular convex cone. It is clear
that $r = m (p + 1), d = (m - 1)$ and $\mV_0$ is the origin of the
convex cone $\mC_0$; hence $\mV_0 \subset \mC_0$. This testing problem
is treated in considerable detail by \cite{pilla2:05}. 

\section{Hypothesis Testing Under Convex Cone Alternatives for
Correlated Data}    
\label{sec-gcone}
In this section, we first derive a statistic for the general testing
problem (\ref{eq:hyp0}) using the decomposition of $\mbg \in
\Re^r$ and next define a new co-ordinate system to transform the null
space $\mV$. Lastly, we derive the asymptotic distribution of the our
test statistic under the model hypothesis.

\subsection{\em Generalized quasi-score statistic for correlated data
and canonical formulation of the testing problem}
\label{sec-can}
Define the {\em generalized quasi-score} (GQS) statistic as
\bea
  \label{eq:sn}
  \S_N := \mQ_N(\overline{\mbg}) - \mQ_N(\widetilde{\mbg}) 
\eea
\noindent for testing the hypothesis (\ref{eq:hyp0}), where
$\overline{\mbg}$ and $\wtmbg$ are defined in Section
\ref{sec-asymp}. 

It is more convenient to define a co-ordinate system to transform the
null space $\mV$. If an appropriate transformation is found, we can
reduce the general problem to a standardized form involving
projections of independently and identically distributed 
standard Gaussian random variables as described in the next section.

\begin{figure}[htp]
\centerline{\scalebox{0.8}{\includegraphics{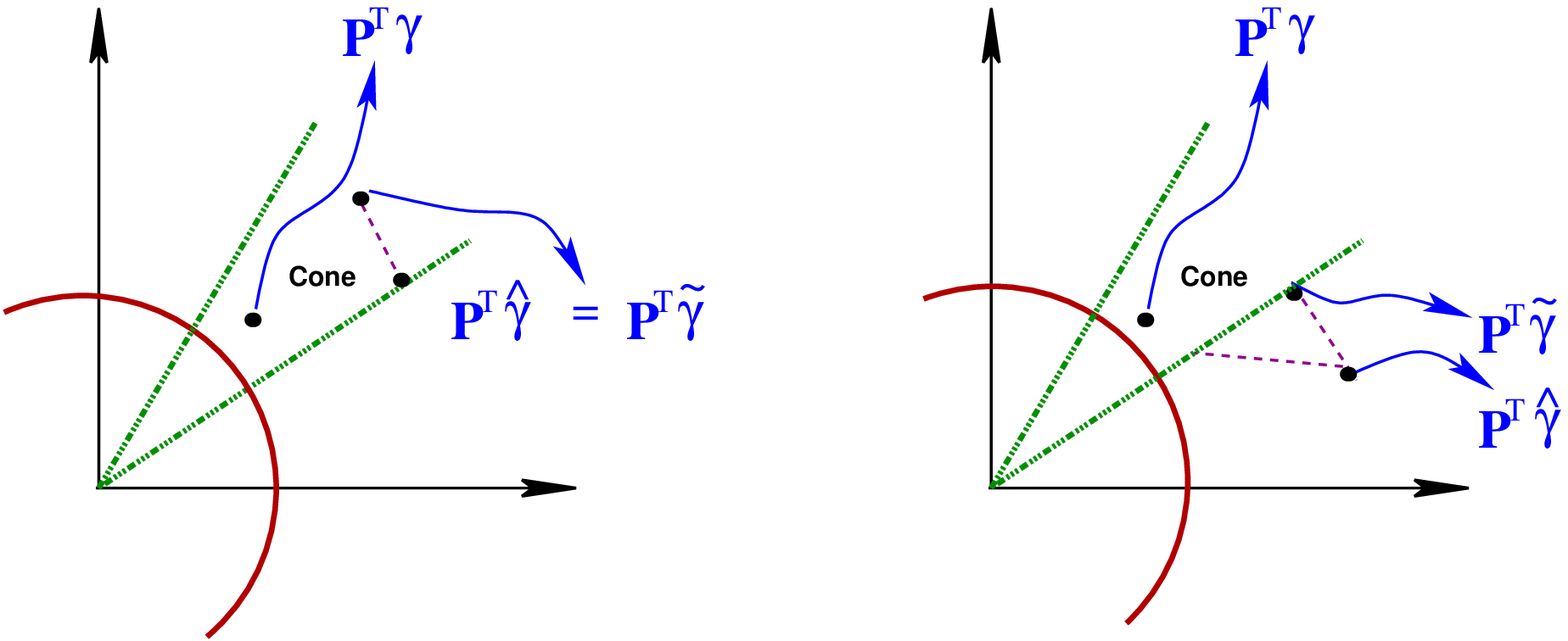}}}
\caption{Projection of $\bP^T \wmbg$ onto the cone $\mC_1$ for
$d = 2$, the number of constraints imposed under $\mH_0^{\, 2}$.} 
\label{fig:cone1}
\end{figure}

Let $\bP$ be a basis matrix for the space $\mV^{\perp}$ whose columns
correspond to the constraints imposed by $\mH_0$. For example, if $d$
constraints are imposed on $\mG$ under $\mH_0$, then $\bP$ is an $(r
\times d)$-dimensional matrix. The choice of the matrix $\bP$ is
problem dependent as shown next.
\begin{lemma}
\label{lemm5}
The hypothesis testing problem (\ref{eq:hyp0}) can equivalently be
represented in terms of the canonical space as testing for 
\bea
 \label{eq:nho}
  \mH_0^{\, 2}\!: \bP^T \mbg = {\bf 0} \quad \mbox{against} \quad
  \mH_1^{\, 2}\!: \bP^T \mbg \in \mC_1 := \bP^T \mC .
\eea
\end{lemma}
\begin{proposition} 
\label{thm2}
{\em Under assumptions \ref{as0}--\ref{as7} and when $\mH_0^{\, 2}\!:
\bP^T \mbg = {\bf 0}$ holds, $\sqrt{N} \, \left(\bP^T \wmbg - \bP^T
\mbg \right) \rightsquigarrow \bZ^{\star} \sim N_{d} [{\b
0}, \Om(\mbg_{0})]$ as $N \rightarrow \infty$, where $\Om(\mbg_0) :=
\bP^T \, \bJ^{-1}(\mbg_{0}) \, \bP$. } 
\end{proposition} 
\begin{example} (Order-restricted testing with three treatment groups). 
\label{eg-pmatrix} In the case of order-restricted testing with three
groups, $\mbg = (\mu_1, \mu_2, \mu_3, \beta_1, \beta_2, \beta_3)^{\,
T}$ and $\mV$ consists of vectors of the form $(\mu, \mu, \mu,
\beta_1, \beta_2, \beta_3)^{\, T}$ which has dimension 4. A basis
matrix for $\mV^{\perp}$ is  
\[
  \bP^T = \begin{pmatrix} 1 \; \; & -1 \;
  \; & 0 \; \; & 0 \; & 0 \; \; & 0 \\ 1 \; \; & 1 \; \; & -2
  \; \; & 0 \; \; & 0 \; \; & 0 \end{pmatrix}
\]
so that
\bea
  \label{eq:bpT}
  \bP^T \mbg = \left( \begin{array}{cc}
  \mu_{_1} - \mu_{_2} \\
  \mu_{_1} + \mu_{_2} - 2 \mu_{_3} \\
  \end{array} \right).
\eea
Therefore, the rows of the matrix $\bP$ span the space
$\mV^{\perp}$. Fig.\ \ref{fig:cone1} demonstrates that (i) if $\bP^T
\wmbg$ lies in the interior of the convex cone $\mC_1$, then $\bP^T
\wmbg = \bP^T \wtmbg$ and (ii) if $\bP^T \wmbg$ lies outside $\mC_1$,
then $\bP^T \wtmbg$ is a projection of $\bP^T \wmbg$ onto $\mC_1$ 
(the orthogonal projection if $\Om(\mbg_0) = \bI_d$). 
\end{example}

\begin{remark} Every $\mbg \in \mG \subset \Re^r$ admits a unique {\em
orthogonal decomposition} of the form $\mbg = \mbg_1 + \mbg_2$ such
that $\mbg_1 \in \mV$ and $\mbg_2 \in \mV^{\perp}$ and $\mC
\subset \mV^{\perp}$. It is clear that $\bP^T \mbg = \bP^T \mbg_2$;
therefore, $\bP^T \mbg \in \bP^T \mC$ since $\bP^T \mV = {\bf 0}$ under
$\mH_0^{\, 2}$.
\end{remark}

Let $\mP_{\mK} \, \bZ$ be the projection of $\bZ$ onto $\mK$
and let $\|\cdot\|$ denote the vector norm.

\subsection{\em Asymptotic equivalence between $\S_N$ and $\|\mP_{\mK} \, 
\bZ\|^2$}  
\label{sec-aequiv}

Owing to Lemma \ref{lemm5}, the hypothesis testing problem
(\ref{eq:hyp0}) is equivalent to that of (\ref{eq:nho}). It will be
established that finding the limiting distribution of $\S_N$ for the
model hypothesis, to appropriate statistical order, is equivalent to
finding the limiting distribution of a length of a certain projection
onto $\mK$.

We present the two main theorems of this section. The first is based
on the quadratic approximation to the inference function $\mQ_N(\mbg)$
in an $N^{-1/2}$-neighborhood of $\mbg_0$ and the second is based on
the transformed null space.
\begin{theorem}
\label{lemm6}
Under assumptions \ref{as0}--\ref{as7}, $\mQ_N(\mbg) -
\mQ_N(\wmbg) =  N \; (\mbg - \wmbg)^{\, T} \; \bJ(\mbg_{0}) \; (\mbg -
  \wmbg) + o_P(1)$ as $N \to \infty$.   
\end{theorem}
{\sc Proof.} Let $\wmbx := \sqrt{N} \, (\wmbg - \mbg_{0})$
so that $\wmbg = (\mbg_{0} + N^{-1/2} \; \wmbx)$. From the quadratic
approximation of the QIF [Theorem 5 of \cite{pilla:05}] the following
result holds
\bea
  \label{qexp}
  \mQ_N \left(\mbg_{0} + N^{-1/2} \, \mbx \right) =
  \mQ_N(\mbg_{0}) + 2 \ip{\mbx}{\bW_N} + \mbx^T \;
  \bJ(\mbg_{0}) \; \mbx + o_P(1),
\eea
where $\mbx \in \Re^r$ is a fixed vector, $\ip{}{}$ is the inner
product and $\bW_N = (2 \sqrt{N})^{-1} \; \nabla
\mQ_N(\mbg_{0})$. Equivalently, 
\bea
  \label{Qxi2}
  \mQ_N(\mbg) - \mQ_N(\mbg_0) &=& 2 \ip{\mbx}{\bW_N} + \mbx^T \;
  \bJ(\mbg_{0}) \; \mbx + o_P(1).
\eea
The minimizer $\mbx^{\star}$ of the quadratic approximation in
(\ref{qexp}) is given by $\mbx^{\star} = - \bJ^{-1}(\mbg_0)
\bW_N$. From Theorem \ref{qnthm}, $\bW_N$ has a limiting distribution
and hence it follows that $\mbx^{\star}$ lies in the ball of radius
${\rm r}_N$ with probability converging to 1. This fact, combined with
the uniformity of the error term in (\ref{qexp}) yields the next
result. If $\wmbx$ is the minimizer of the quadratic approximation in
(\ref{qexp}), then the QIF estimator becomes $\wmbg = (\mbg_0 +
N^{-1/2} \, \wmbx)$. Equivalently, $\wmbx = -\bJ^{-1}(\mbg_0) \, \bW_N
+ o_P(1)$. Therefore, it follows that
\bea
  \label{eq:qnalt}
  \mQ_N(\wmbg) - \mQ_N(\mbg_0) &=& 2 \ip{\wmbx}{\bW_N} + \wmbx^T
  \; \bJ(\mbg_{0}) \; \wmbx + o_P(1).
\eea
Since $\bW_N = -\bJ(\mbg_{0}) \; \wmbx + o_P(1)$, equations
(\ref{Qxi2}) and (\ref{eq:qnalt}) simplify to $\mQ_N(\mbg) -
\mQ_N(\wmbg) = (\mbx - \wmbx)^T \; \bJ(\mbg_{0})  \; (\mbx - \wmbx) +
o_P(1) = N \, (\mbg - \wmbg)^T \; \bJ(\mbg_{0}) \; (\mbg - \wmbg) +
o_P(1)$. \hfill \rule{2mm}{2mm} 

In the next theorem, we establish the relation between the GQS
statistic for correlated data and the squared length of projection of
the standard Gaussian vector $\bZ$ onto $\mK$ for independent
data. Consequently, we can derive a result for the asymptotic null
distribution of $\S_N$ by employing the seminal work of
\cite{take:97}.
\begin{theorem}
\label{thm-proj}
Under assumptions \ref{as0}--\ref{as8}, the GQS statistic $\S_N$
for testing the hypothesis under a general convex cone alternative 
(\ref{eq:hyp0}) is asymptotically equivalent to the statistic
$\|\mP_{\mK} \, \bZ\|^2$ for testing the hypothesis (\ref{eq:nho}),
where $\bZ \sim N_r({\bf 0}, \bI_r)$. That is, 
\bea
  \label{eq:local}
  \S_N \rightsquigarrow \|\mP{_\mK} \,
  \bZ\|^2 \quad \mbox{as} \quad N \to \infty.
\eea
\end{theorem} 
{\sc Proof.}  The QIF estimators obtained by minimizing
$\mQ_N(\mbg)$ under the spaces $\mV, (\mV \oplus \mC)$ and $\mG
\subset \Re^r$ are ordered as $\mQ_N(\overline{\mbg}) \geq
\mQ_N(\wtmbg) \geq \mQ_N(\wmbg)$. From Theorem \ref{lemm6}, in an
$N^{-1/2}$-neighborhood of the true parameter vector $\mbg_{0}$,
$\mQ_N(\mbg) - \mQ_N(\wmbg) = N \; (\mbg - \wmbg)^{\, T} \;
\bJ(\mbg_{0}) \; (\mbg - \wmbg) + o_P(1)$. From Proposition  
\ref{thm2}, it suffices to consider the transformed hypothesis
(\ref{eq:nho}). There exists an $(r \times r)$-dimensional matrix
$\bL$ such that $\bL^T \bL = \bJ(\mbg_0)$; for example, Cholesky
factorization of $\bJ(\mbg_0)$.

It is more convenient to consider $\mbr$-parametrization under the 
transformed null space $\bL \, \mV =: \mV^{ \star}$, where $\mV$ is 
the null space under the $\mbg$-parametrization. Let $\mbr =
\sqrt{N} \, \bL \, \mbg$ so that $\wmbr = \sqrt{N} \, \bL \,
\wmbg, \, \overline{\mbr} = \sqrt{N} \, \bL \, \overline{\mbg}$ and
$\wtmbr = \sqrt{N} \, \bL \, \wtmbg$. By Theorem \ref{thm12}, 
\bea
  \label{eq:mbr} 
  (\wmbr - \mbr_{0}) \rightsquigarrow \bZ
  \sim N_r({\bf 0}, \bI_r) \quad \mbox{as} \quad N \to \infty.
\eea
Furthermore, Theorem \ref{lemm6} yields 
\bea
  \label{eq:nqn} 
  \mQ_N (\mbg) &=& \mQ_N(\wmbg) + \| \mbr - \wmbr\|^2 +
  o_P(1). 
\eea

A given $\mbr \in \Re^r$ admits an orthogonal decomposition as $\mbr =
\mbr_1 + \mbr_2$ such that $\mbr_1 \in \mV^{\star}$ and 
$\mbr_2 \in (\mV^{\star})^{\perp}$. However, orthogonality is not
preserved by $\bL$; hence $(\mV^{\star})^{\perp} \ne \bL \,
\mV^{\perp}$. The hypothesis (\ref{eq:nho}) can be re-expressed as
\bea
  \label{eq:tho}
  \mH_0^{\, 3}\!: \mbr_1 \in \mV^{ \star}, \, \mbr_2 = {\bf 0} \quad
   \mbox{against} \quad \mH_1^{\, 3}\!: \mbr_1 \in \mV^{\star}, \,
  \mbr_2 \in \mK \quad \mbox{and} \quad \mbr_2 \neq {\bf 0},
\eea
where $\mK = [\bL \, (\mC \oplus \mV) \} \cap
\{(\mV^{\star})^{\perp}]$ is a $d$-dimensional cone in $\Re^r$ since
$\mbox{dim}[(\mV^{\star})^{\perp}] = d$.  

Similarly, the estimator $\wmbr = \sqrt{N} \, \bL \, \wmbg$ admits an
orthogonal decomposition as $\wmbr_1 + \wmbr_2$ such that $\wmbr_1 \in
\mV^{\star}$ and $\wmbr_2 \in (\mV^{\star})^{\perp}$. By
orthogonality, equation (\ref{eq:nqn}) can be expressed as 
\bea
  \label{qn1}
  \mQ_N \left[ N^{-1/2} \, \bL^{-1} \, (\mbr_1 + \mbr_2)
  \right] - \mQ_N \left(\wmbg\right) = \|\mbr_1 - \wmbr_1\|^2 +
  \|\mbr_2 - \wmbr_2\|^2 + o_P(1). 
\eea 
Under $\mH_0^{\, 3}$, (\ref{qn1}) simplifies to
\bea
  \label{eq:min}
  \mQ_N \left[ N^{-1/2} \, \bL^{-1} \, (\mbr_1 + \mbr_2)
  \right] - \mQ_N(\wmbg) = \| \mbr_1 - \wmbr_1\|^2 + \| \wmbr_2 \|^2
  + o_P(1).  
\eea
Under $\mH_{0}^{\, 3}$, $\mbr_1 \in \mV^{\star}$ and $\mbr_2 = {\bf
0}$, whereas under $\mH_{1}^{\, 3}$, $\mbr_1 \in
\mV^{\star}$ and $\mbr_2 \in \mK$.  

At first, consider minimizing over $\mH_{0}^{\, 3}\!: \mbr_2 = 0$. The
right-hand side of (\ref{eq:min}) is uniquely minimized at $\mbr_1 =
\wmbr_1$ and $\mbr_2 = 0$. By definition, the left hand side is 
minimized at $\mbr_1 = \overline{\mbr}_1$. By uniformity of the error
term and uniqueness of the minimum, it follows that $\overline{\mbr}_1
= \wmbr_1 + o_P(1)$. 
\begin{figure}[htb]
\centerline{\scalebox{0.8}{\includegraphics{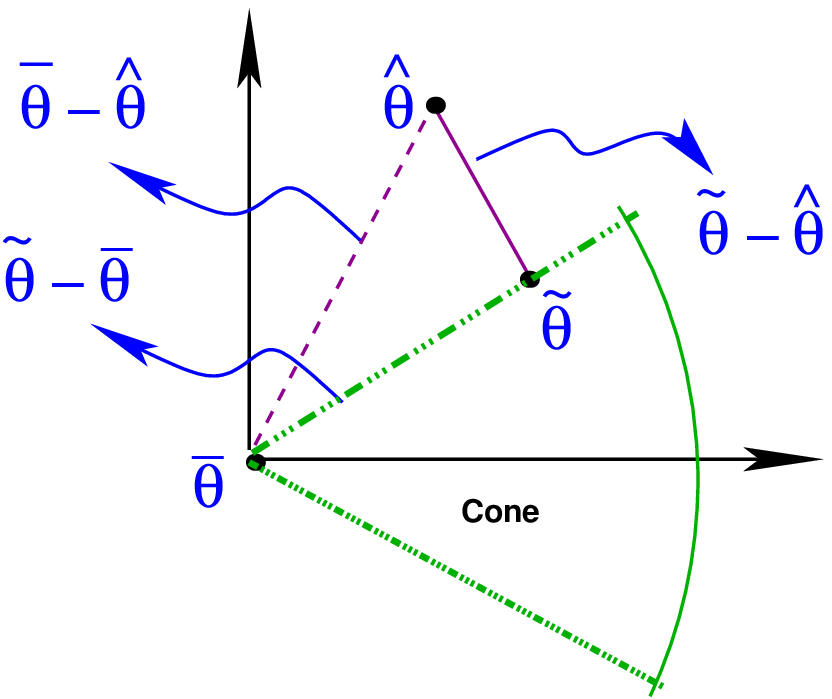}}}
\caption{Projections of estimators of $\mbr$ onto the cone 
$\mK$. Asymptotically, the relation $\wtmbr = \mP_{\mK} \, \wmbr + 
o_P(1)$ holds.}  
\label{fig:cone2}
\end{figure}
Under $\mH_1^{\, 3}$, the right-hand side of (\ref{qn1}) is minimized
at $\mbr_1 = \wmbr_1$ and $\mbr_2 = \mP_{\mK} \, \wmbr_2$. Therefore,
the left hand side of (\ref{qn1}) is minimized at $\wtmbr_1 =
\wmbr_1 + o_P(1)$ and $\wtmbr_2 = \mP_{\mK} \, \wmbr_2 +
o_P(1)$. In effect, minimizing $\mQ_N(\cdot)$ in (\ref{qn1}) under
$\mH_0^{\, 3}$ and $\mH_1^{\, 3}$ yield respectively, $\mQ_N
(\overline{\mbg}) - \mQ_N(\wmbg) = \|\wmbr_2\|^2 + o_P(1)$ and $\mQ_N
(\wtmbg) - \mQ_N(\wmbg) = \|\wtmbr_2 - \wmbr_2\|^2 + o_P(1)$ since $\|
\overline{\mbr}_1 - \wmbr_1\|^2 = o_P(1)$ and $\| \wtmbr_1 -
\wmbr_1\|^2 = o_P(1)$. Consequently, $\S_N = \left[
\mQ_N(\overline{\mbg}) - \mQ_N(\wtmbg) \right] = \| \wmbr_2\|^2 -
\|\wtmbr_2 - \wmbr_2\|^2 + o_P(1) = \| \wmbr_2\|^2 - \|\mP_{\mK} \,
\wmbr_2 - \wmbr_2\|^2 + o_P(1) = \| \mP_{\mK} \, \wmbr_2 \|^2 +
o_P(1)$ since $\left(\mP_{\mK} \, \wmbr_2 - \wmbr_2 \right)$ and
$\mP_{\mK} \, \wmbr_2$ are orthogonal as shown in Fig.\
\ref{fig:cone2}. Furthermore, since $\mP_{\mK} \, \wmbr_1 = 0$ under 
$\mH_0^{\, 3}$ given that $\wmbr \in  \mV^{\star}$ and $\mK \subset
(\mV^{\star})^{\perp}$, it follows that $\S_N = \| \mP_{\mK} \, \wmbr
\|^2 + o_P(1)$. Moreover, $\mP_{\mK} \, \mbr_0 = 0$ yielding
$\mP_{\mK} \, \wmbr = \mP_{\mK} \, (\wmbr - \mbr_0)$. The result
follows from (\ref{eq:mbr}). \hfill \rule{2mm}{2mm} 

\subsection{Asymptotic null distribution of the generalized
quasi-score test}  
\label{sec-adist} 
We derive the asymptotic distribution of the GQS statistic $\S_N$ when
$\mH_0^{\, 2}$ holds. Until now, we established that (1) the testing
problems (\ref{eq:hyp0}) and (\ref{eq:nho}) are equivalent (Lemma
\ref{lemm5}, Section \ref{sec-can}) and (2) there exists an asymptotic
relation between the statistic $\S_N$ based on correlated data for the
testing problem (\ref{eq:nho}) and the statistic $\|\mP_{\mK} \,
\bZ\|^2$ based on independent data for the testing problem 
(\ref{eq:hyp1}) (Theorem \ref{thm-proj}, Section
\ref{sec-aequiv}). These main results, in conjunction with Theorem 2.1
of \cite{take:97}, yield the asymptotic null distribution of $\S_N$
for the testing problem (\ref{eq:nho}). The weights of this asymptotic
null distribution are mixed volumes of $\mK$ and its {\em polar} or
{\em dual cone} $\mK^{0}$ \citep{webster:94}.

Let $\mS^{(d - 1)}$ be the closed $(d - 1)$-dimensional unit sphere in
$\Re^r$, $\mM = \mK \cap \mS^{(d - 1)}$ be the $(d - 1)$-dimensional
convex manifold and $\mM^{0} = \mK^{0} \cap \mS^{(d - 1)}$. Let
$\vartheta_{d - k, k}(\mM, \mM^{0})$ be the {\em mixed volumes} of
$\mM$ and $\mM^{0}$ for $k = 0, \ldots, d$.
\begin{theorem}
\label{them-adist}
Under assumptions \ref{as0}--\ref{as8} and when $\mH_0\!: \mbg \in
\mV$ holds, the asymptotic distribution of the GQS statistic $\S_N$
for any $c > 0$ is  
\bea
 \label{eq:adist} 
 \underset{N \rightarrow \infty}{\lim} \mbP \left(
 \S_N \leq c \right) &=& \mbP \left(\| \mP_{\mK} \, \bZ \|^2 \leq
 c \right) \nonumber \\ &=& \sum_{k = 0}^{d} \binom{d}{k} \,
 \frac{\vartheta_{d - k, k} \left(\mM, \mM^{0} \right)}{\omega_{_k} \;
 \omega_{_{d - k}}} \; \mbP \left(\chi^2_{d - k} \leq c \right),
\eea
where $d := \mbox{dim}(\mK)$, $\omega_{_{k - 1}} = 2 \,
\pi^{k/2}/\Gamma(k/2)$ is the volume of the unit sphere $\mS^{(k -
1)}$ embedded in $\Re^r$ and $\chi^2_k$ is the chi-squared
distribution with $k$ degrees of freedom. The $\chi^2_k$ for $k = 0$
is simply a point mass at the origin so that $\mbP(\chi^2_0 \leq c) =
1$ for all $c > 0$.   
\end{theorem}

The right-hand side of (\ref{eq:adist}) is a weighted mean of several
tail probabilities of chi-squared distributions; hence, it is often
referred to as {\em chi-bar-squared} distribution and denoted by
$\overline{\chi}^2$ \citep{shapiro:88}. In general, it is very
difficult to derive explicit expressions for the weights for the
asymptotic null distribution in (\ref{eq:adist}). In certain special
cases of $\mK$, weights are known explicitly or can be evaluated
numerically. For instance, for polyhedral cones (i.e., the cones
defined by a finite number of linear constraints) one can calculate
the weights. For the general case of non-polyhedral cones, Section 3.5
of \cite{silva:04} provide a simulation-based approach.

\cite{take:97} clarify the geometric meaning of the weights
when the boundary of the cone is smooth or piecewise smooth. However,
as the following example demonstrates, in the case of order restricted
testing problem, one does not have a smooth cone or a smooth manifold
and hence a more general approach to calculate the weights is
warranted which is derived in Section \ref{sec-wts}. \\
\begin{example} (Asymptotic null distribution of $\S_N$ for three
treatment groups). \label{eg-angle} Consider the problem of order
restricted testing with three treatment groups. Under $H_0^o\!:
\mu_{_1} = \mu_{_2} = \mu_{_3}$, the asymptotic distribution in
(\ref{eq:adist}) has an explicit expression. The convex cone $\mK$ is
the region between the two vectors defining the constraints in
(\ref{eq:bpT}). Let $\phi$ be the angle of the cone $\mK \subset
\Re^2$ at the vertex. Derivation of $\phi$ will be presented in Example
\ref{eg-coneN} in Section \ref{sec-prep}. Note that the angles (in
radians) of $\mK$ and $\mK^0$ at their vertices sum to $\pi$. Divide
the plane into the following four regions: 
\bee
\vspace{-0.15in}
\item The cone $\mK$ such that $\bZ \in \mK$ with a probability of
$\phi/(2 \, \pi)$ and conditional on $\bZ \in \mK$, it follows that
$\mP_{\mK} \, \bZ = \bZ$ with $\|\mP_{\mK} \, \bZ\|^2 \sim \chi_2^2$. 
\vspace{-0.15in}
\item The dual cone $\mK^{0}$ such that $\mP_{\mK} \, \bZ = 
0$ for all $\bZ \in \mK^{0}$ with a probability of $[1/2 - \phi/(2
\pi)]$ and conditional on $\bZ \in \mK^{0}$, it follows that
$\|\mP_{\mK} \, \bZ\|^2 \sim \chi_0^2$. 
\vspace{-0.15in}
\item The two regions $\mK^{\dag}$ and $\mK^{\star}$, 
where $\mP_{\mK} \, \bZ$ is a multiple of one of the two vectors
defining the constraints in (\ref{eq:bpT}); conditional on $\bZ \in
\mK^{\dag}$ or $\bZ \in \mK^{\star}$ with a total probability of
$1/2$, it follows that $\|\mP_{\mK} \, \bZ\|^2 \sim \chi^2_1$. 
\vspace{-0.15in}
\eee
For $d = 2$, the asymptotic distribution in (\ref{eq:adist})
simplifies to a $\overline{\chi}^2$ distribution:
\bea
  \label{eq:SN0}
  \S_N \rightsquigarrow \|\mP_{\mK} \, \bZ\|^2 \sim
  \left(\frac{1}{2} - \frac{\phi}{2 \, \pi} \right) \chi_0^2 +
  \frac{1}{2} \chi_1^2 + \frac{\phi}{2 \, \pi} \chi_2^2 \quad
  \mbox{as} \quad N \to \infty. 
\eea
\end{example}

\begin{remark}
The right hand side expression of (\ref{eq:SN0}) also occurs in the
context of the asymptotic null distribution of the LRT statistic for
testing for two-component mixture model
\citep{lindsay:95,lin2:97,pilla:03}. 
\end{remark}

\begin{example} (Asymptotic null distribution of $\S_N$ under order 
restricted alternatives). We consider the order restricted testing
problem discussed in Section \ref{sec-order}. For convenience, we
reorder the elements of the parameter vector as $\mbg = (\mbm^T,
\mbt^T)^T$, with corresponding permutations of the rows and columns of
$\bJ(\mbg_0)$ [see Lemma \ref{lemm2} for the definition of
$\bJ(\mbg_0)$]. We partition $\bJ(\mbg_0)$ as
\bes
  \bJ(\mbg_{0}) := \begin{pmatrix} \bJ_{\mbm \mbm} \; \; & \bJ_{\mbm
  \mbt} \cr \bJ_{\mbt \mbm} \; \;  & \bJ_{\mbt \mbt}
\end{pmatrix}.
\ees
Let $\bJ^{\mbm \mbm}$ be the appropriate submatrix of
$\bJ^{-1}(\mbg_0)$. From the formula for an inverse of a partitioned
matrix, it follows that $\bJ^{\mbm \mbm} = \left(\bJ_{\mbm \mbm} -
\bJ_{\mbm \mbt} \, \bJ_{\mbt \mbt}^{-1} \, \bJ_{\mbt \mbm}
\right)^{-1}$. Since subjects in different groups are independent, the
variance matrix $\bJ^{\mbm \mbm}$ is diagonal. Let $\bZ^{\dag} \sim
N_m({\bf 0}, \bQ)$, where $\bQ$ is a pre-specified diagonal variance
matrix.

In the order restricted testing problem, the asymptotic null
distribution has an explicit expression. Under $\mH_0^{\, o}$, for any
$c > 0$, the result (\ref{eq:adist}) reduces to $\underset{N
\rightarrow \infty}{\lim} \mbP \left( \S_N \leq c \right) = \sum_{k =
0}^{d} \; p(d - k + 1, d; \bQ) \; \mbP \left(\chi^2_{d - k} \leq c
\right)$, where $\bQ = \bJ^{\mbm \mbm}$ and $p(k, d; \bQ)$ is the {\em
level probability} that the projection of $\bZ^{\dag}$ onto $\mC_0$
with a weight vector $\bQ$ consists of exactly $k$ distinct points
[Section 2.4, \cite{robert:88}]. The unknown weight vector is replaced
with $\widehat{\bQ}$. This problem is developed and treated in
considerable detail by \cite{pilla2:05}. The weights $p(k, d; \bQ)$
are also referred to as $\overline{\chi}^2$ weights
\citep{robert:88}. Such weights appear in the null asymptotic or exact
distribution of several test statistics when there are inequality
constraints on parameters. In certain cases, exact expressions for
these weights are available and in other cases, one may obtain
approximations or bounds \citep{silva:04}. 
\end{example}

\section{Asymptotic Null Distribution of $\S_N$: The Volume-of-Tube
Formula} 
\label{sec-wts}
In this section, we derive explicit expressions for the weights in the
asymptotic null distribution of $\S_N$ in (\ref{eq:adist}) by 
representing $\mM$ and $\mK$ in parametric form and in turn using the 
Hotelling-Weyl-Naiman volume-of-tube formula.

\subsection{Parametric representation of $\mM$ and $\mK$}
\label{sec-prep} From a geometrical perspective, the rows of the $(r
\times d)$-dimensional matrix $\bP$ span the space $\mV^{\perp}$. From
the orthogonal decomposition in Section \ref{sec-gcone}, it follows
that $\mbg = \mbg_1 + \mbg_2 = \mbg_1 + \bP \mbn$ such that $\mbg_1
\in \mV$ and for some $\mbn \in \Re^d$. Let $\mN := \{\mbn: \bP \mbn
\in \mC\}$. Given that $\mbg \in \mC = \{\bP \mbn\!: \mbn \in \mN\}$,
the hypothesis (\ref{eq:hyp0}) can be expressed as
\bea
  \label{hypn}  
  \mH_0^{\, 4}\!: \mbn = {\bf 0} \quad \mbox{against} \quad \mH_1^{\,
  4}\!: \mbn \in \mN.  
\eea

In order to represent $\mK$ and $\mM$ in a parametric form, we require
the following result.
\begin{proposition}  
\label{prop2}
The matrix $\bP^{\star} := (\bL^{-1})^T \, \bP$ forms a basis for
$(\mV^{\star})^{\perp}$, where $\mV^{\star} = \bL \, \mV$.
\end{proposition}  
{\sc Proof.} Let $\b y \in \mV^{\star}$ so that $\b y = \bL \b z$ for
some $\b z \in \mV$. Let $\b x \in C(\bP^{\star})$, the column space
of the matrix $\bP^{\star}$, so that $\b x \in (\bL^{-1})^T \bP \mbn$
for some vector $\mbn = (u_{_1}, \ldots, u_{_d})^T \in \Re^d$. It
follows that, $\ip{\b x}{\b y} = (\mbn^T \bP^T \bL^{-1}) (\bL \b z) =
\ip{\bP \mbn}{\b z} = 0$ since $\b z \in \mV$ and $\bP \mbn \in
\mV^{\perp}$. Therefore, the space spanned by $\bP^{\star}$ is a
subspace of $(\mV^{\star})^{\perp}$. However, this subspace and
$(\mV^{\star})^{\perp}$ have the same dimension $d$, therefore they
must be equal. \hfill \rule{2mm}{2mm} 

We return to the orthogonal decomposition of $\mbr = \mbr_1 +
\mbr_2$ such that $\mbr_1 \in \mV^{\star}$ and $\mbr_2 \in
(\mV^{\star})^{\perp}$. Proposition \ref{prop2} ensures the following 
representation: $\mbr_2 = \mP_{(\mV^{\star})^{\perp}} \, \mbr =
\mP_{\bP^{\star}} \, \mbr = \bP^{\star} \left[ \left(\bP^{\star}
\right)^T \bP^{\star} \right]^{-1} \, \left(\bP^{\star} \right)^T \mbr
= \bP^{\star} \, \left[\bP^T \, \bJ(\mbg_0) \, \bP \right]^{-1} \,
\bP^T \, \bL^{-1} \mbr$ since $\bL^T \bL = \bJ(\mbg_0)$. Therefore,
\bea
 \mbr_2 &=& \sqrt{N} \, \bP^{\star} \, \left[\bP^T \, \bJ(\mbg_0) \,
  \bP \right]^{-1} \, \bP^T (\mbg_1 + \bP \mbn) \nonumber \\
  \label{eq:theta2} 
  &=& \sqrt{N} \, \bH \, \mbn \quad \mbox{since} \quad \bP^T \mbg_1 =
  {\bf 0}, 
\eea
where 
\bea
  \label{eq:hmat}
  \bH := \bP^{\star} \, [\bP^T \, \bJ(\mbg_0) \, \bP]^{-1} \, \bP^T
  \, \bP. 
\eea
Next, the $d$-dimensional cone $\mK = \{\bL \, (\mC
\oplus \mV) \} \cap \{(\mV^{\star})^{\perp}\}$ and the $(d -
1)$-dimensional manifold $\mM = \mK \cap \mS^{(d - 1)}$ are
represented in the parametric form by considering a vector function
$\mT(\mbn)\!: \mN \subset \Re^{d} \rightarrow \mM \subset \Re^{(d -
1)}$, where $\mbn \in \mN$. That is,
\bea
  \label{eq:Tn}
  \mT(\mbn) := \{ \bH \, \mbn\!: \mbn \in \mN \subset \Re^d \;
  \mbox{and} \; \|\bH \, \mbn\| = 1\}, 
\eea
where $\bH$ is an $(r \times d)$ matrix defined in (\ref{eq:hmat}). In
effect, we can redefine $\mM := \{\mT(\mbn) \in \mS^{(d - 1)}\!: \mbn
\in \mN \subset \Re^d\}$ and $\mK := \{\bH \, \mbn\!: \mbn \in \mN
\subset \Re^d\}$.  \\  
\begin{example} (Explicit Expressions for the convex cone $\mN \subset
\Re^2$ and the angle $\phi$ of $\mK$).
\label{eg-coneN} We return to the problem of testing under order
restricted hypothesis for correlated data. If $\mH_1^{\, o}\!:
\mu_{_1} > \mu_{_2} > \mu_{_3}$, then the choice of $\bP$ is given in
Example \ref{eg-pmatrix}. Therefore, $\mbg_2 = \bP \mbn$ corresponds
to $\mu_{_1} = (u_{_0} + u_{_1} + u_{_2}), \mu_{_2} = (u_{_0} -
u_{_1} + u_{_2})$ and $\mu_{_3} = (u_{_0} - 2 u_{_2})$, where 
$u_{_0} = (\mu_{_1} + \mu_{_2} + \mu_{_3})/3$. The constraints
$\mu_{_1} > \mu_{_2}$ and $\mu_{_2} > \mu_{_3}$ under $\mH_1^{\, o}$
yield respectively, $u_{_1} > 0$ and $u_{_2} >
u_{_1}/3$. Consequently, the convex cone $\mN = \{\mbn\!:
u_{_1} > 0, u_{_2} > u_{_1}/3\}$.  It is clear that
$\mbn$ lies in the cone bounded by the vectors ${\bf v}_{1} = (1,
1/3)^T$ and ${\bf v}_{2} = (0, 1)^T$. The cone $\mK$ is then bounded
by the $\mbr_2$-component of $\bL \, \bP \, {\bf v}_1$ and $\bL \, \bP
\, {\bf v}_2$. Hence  
\bes
  \cos(\phi) = \ip{\frac{\bH {\bf v}_{1}}{\|\bH {\bf v}_{1}\|}}{
  \frac{\bH {\bf v}_{2}}{\|\bH {\bf v}_{2}\|}}
\ees
yields an explicit expression for the angle $\phi$ defined in Example
\ref{eg-angle}.    
\end{example}

\subsection{\em Asymptotic null distribution of $\S_N$ in terms of the
geometry of $\mM$}
As a first step, we establish the connection between the distribution
of a squared length of projection of $\bZ$ onto $\mK$ and the
volume-of-tube problem. We take a different approach from
\cite{lin:97} in order to cast the problem in the general framework of
this article. 

The geodesic (or angular) distance between two points on any manifold
is defined as the shortest measured distance between the points {\em
within} the manifold itself. Let $\fT(\varrho, \mM)$ or $\fT(\phi,
\mM)$ be the spherical {\em tube} around the topological $(d -
1)$-dimensional manifold $\mM$ of Euclidean radius $\varrho$ or
geodesic radius $\phi$ embedded in $\mS^{(d - 1)}$, where $\varrho =
\sqrt{2 [1 - \cos(\phi)]}$, $\mbn \in \mN$ and ${\rm dim}(\mN) =
d$. Since $\mS^{(d - 1)}$ is also a manifold and $\fT \subset \mS^{(d
- 1)}$, the geodesic distance between two points on $\fT$ is the
length of the segment of the great circle (arc) connecting the two
points. We view each ray $\{\zeta \, \mbe\!: \zeta \geq 0\}$ as a cone
on which to make a projection, yielding $\wZ_{\mbe}$ that depends on
$\mbe$. We redefine the cone as $\mK := \{\zeta \, \mbe\!: \zeta > 0,
\, \|\mbe\| = 1\}$ to yield $\mP_{\mK} \, \bZ = \underset{\zeta \,
\mbe \, \in \, \mK}{\arg \, \min} \; \| \zeta \, \mbe - \bZ\|^2 =
\underset{\mbe \, \in \, \mK}{\sup} \, \ip{\mbe}{\bZ}_+$, where
$\ip{\cdot}{\cdot}_+$ denotes the positive part of the inner
product. In effect, we have 
\bea
 \label{pksup}
 \| \mP_{\mK} \, \bZ \|^2 = \left[\sup_{\mbe \, \in \, \mK} \,
 \ip{\mbe}{\bZ}_+ \right]^2.
\eea

In order to reduce the problem to that of a uniform process, we
condition on $\|\bZ\|^2$ and integrate over the conditional
distribution. Consequently, from (\ref{pksup}) we can express 
\bea
  \mbP \left(\| \mP_{\mK} \, \bZ \|^2 \leq c \right) &=& 
  \mbP \left[ \left( \underset{\mbe \, \in \, \mK}{\sup}
  \ip{\mbe}{ \frac{\bZ}{\|\bZ\|} }_+ \right)^2 \leq \frac{c}{\|
  \bZ \|^2} \right] \nonumber \\
  &=& \int_{c}^{\infty} \mbP \left[ \left( \underset{\mbe \, \in \, 
    \mK}{\sup} \ip{\mbe}{\bU}_+ \right)^2 \leq \frac{c}{z} \Bigg| \| 
  \bZ\|^2 = z \right] \; f_r(z) \, dz \nonumber \\
  \label{eq:intp} 
  &=& \int_{c}^{\infty} \mbP \left( \underset{\mbe \, \in \, 
  \mK}{\sup} \ip{\mbe}{\bU}_+ \leq \sqrt{\frac{c}{z}} \,
  \right) \; f_r(z) \, dz,
\eea
where $\bU = (Z_1/\|\bZ\|, \ldots, Z_r/\|\bZ\|)$ is uniformly
distributed on $\mS^{r}$, an $r$-dimensional unit sphere embedded in
$\Re^{(r + 1)}$, and $f_r(z)$ is a $\chi^2$ density with $r$ degrees
of freedom. Therefore, the right-hand side of (\ref{eq:adist}) can be
determined from (\ref{eq:intp}), provided the probability in the
integrand can be found, at least approximately. 

The uniformity property of $\bU$ reduces the problem of finding $\mbP
(\sup_{\mbe} \ip{\mbe}{\bU}_+ \leq \sqrt{c/z})$ to that of determining
the volume of the tube $\fT(\varrho, \mM)$ including the end points
corrections proposed by \cite{naiman:90}. Consequently,
\bea
  \label{eq:vol}
  \mbP \left( \underset{\mbe \, \in \, \mK}{\sup}
  \ip{\mbe}{\bU}_+ \leq \sqrt{\frac{c}{z}} \, \right) =
  \mbP \left[ \bU \in \fT(\phi, \mM) \right] =
  \frac{\vartheta_{\mM}(\phi)}{\omega_{_{r - 1}}}, 
\eea
where $\omega_{_{r - 1}} = 2 \, \pi^{r/2}/\Gamma(r/2)$ is the volume
of $\mS^{r}$ and $\vartheta_{\mM}(\phi)$ is the volume of
$\fT(\varrho, \mM)$. Therefore, (\ref{eq:intp}) and (\ref{eq:vol})
establish a connection between $\|\mP_{\mK} \, \bZ\|^2$ and volume of
the tube $\fT(\phi, \mM)$ around $\mM$ embedded in
$\mS^{r}$. Essentially, we established that the distribution of $\|
\mP_{\mK} \, \bZ \|^2$ can be determined explicitly by finding 
$\vartheta_{\mM}(\phi)$ which equals $\mbox{cos}^{-1}(\sqrt{c/z})$ for
any $0 \leq \phi \leq \pi/2$. 

The asymptotic expansion of the tail probability of the ${\sup}
\ip{\mbe}{\bZ}_+$ can also be obtained using the {\em Euler-Poincar\`e
characteristic} $\mathcal{E}$ method, developed by \cite{adler:81} and
\cite{worsley2:95,worsley1:95}, where the expectation of the
$\mathcal{E}$ of an excursion set is evaluated. \cite{take:02}
establish the equivalence between the tube and Euler characteristic
methods under the assumption that $\mM$ is a manifold with piecewise
smooth boundary.

We motivate the geometric concepts through the order restricted
alternatives for correlated data. We define a {\em corner} to mean a
point where two faces of the boundary of $\mM$ meet.  \\
\begin{example} (Geometry of $\mM$ for the order restricted
alternatives). \label{eg-2case} We assume that $n_i = n \; (i = 1,
\ldots, N)$ and the number of subjects in each group is equal so that
we have a balanced design. In this case, $\omega_{_0} = 2$ and
$\omega_{_1} = 2 \, \pi$. First consider three treatment groups (i.e.,
$m = 3$), then the number of restrictions $d$ equals two corresponding
to $\mu_{_1} < \mu_{_2} < \mu_{_3}$. Therefore, $\mbox{dim}(\mK) = 2$
and $\mM$ is just an arc with two end points. Suppose $m = 4$
corresponding to three constraints, then $\mM$ is a spherical
triangle. The interior corresponds to $\mu_{_1} < \mu_{_2} < \mu_{_3}
< \mu_{_4}$ with three corners $\mu_{_1} = \mu_{_2} = \mu_{_3},
\mu_{_2} = \mu_{_3} = \mu_{_4}, \mu_{_1} = \mu_{_2}$ and $\mu_{_3} =
\mu_{_4}$ and three edges $\mu_{_1} = \mu_{_2}, \mu_{_2} = \mu_{_3},
\mu_{_3} = \mu_{_4}$.    
\end{example}

Example \ref{eg-2case} demonstrates that determining
$\vartheta_{\mM}(\phi)$ for $d \geq 3$ depends on the geometry of
$\mM$. \cite{naiman:90} derived expressions for the volume of a tube
by decomposing the tube into different sections, corresponding to the
main part of the manifold, hemispherical caps along boundaries of the
manifold, circular wedges at the boundaries and so on. Adding up these
terms yields a series involving partial beta functions [Lemma 3.6 of
\cite{naiman:90}]. Substituting these terms into (\ref{eq:intp}) 
yields a series involving partial gamma functions; the first four
terms of which are given in the next theorem whose proof essentially
follows from \cite{pilla:03} and hence is omitted.
\begin{theorem}
\label{th:sig3d}
Under assumptions \ref{as0}--\ref{as8} and when $\mH_0\!: \mbg \in
\mV$ holds, the asymptotic distribution of $\S_N$ for a
$d$-dimensional manifold $\mM$ and for any $c > 0$ is given by 
\bea
   \underset{N \rightarrow \infty}{\lim} \mbP \left(
   \S_N \geq c \right) &=& \mbP \left(\| \mP_{\mK} \, \bZ \|^2 \geq c
   \right) \nonumber \\ 
   &=& \, \frac{\kappa_{_0}}{\omega_{_{d - 1}}} \, \mbP \left(
   \chi_{d}^2 \ge c \right) + \, \frac{\ell_0}{2 \, \omega_{_{d - 2}}}
   \, \mbP \left( \chi_{d - 1}^2 \ge c \right) + \, \frac{(\kappa_{_2}
   + \ell_1 + \upsilon_{_0})}{2 \, \pi \, \omega_{_{d - 3}}} \, \mbP
   \left( \chi_{d - 2}^2 \ge c \right) \nonumber \\
   \label{eq:wts}
   && + \, \frac{(\ell_2 + \upsilon_{_1} + \tau)}{4 \, \pi \,
   \omega_{_{d - 4}}} \, \mbP \left( \chi_{d - 3}^2 \ge c \right)
   + o\left(c^{(d - 5)/2} e^{-c/2} \right) \quad \mbox{as} \quad c \to
   \infty,
\eea
where $\kappa_{_0}$ is the \((d - 1)\)-dimensional volume of the
manifold $\mM$, $\kappa_{_2}$ is the measure of curvature of $\mM$, 
$\ell_0$ is the \((d - 2)\)-dimensional volume of the boundaries of
$\mM$, $\ell_1$ is the measure of rotation of the boundary, $\ell_2$
is the measure of curvature similar to $\kappa_{_2}$, $\upsilon_{_0}$
is the measure of rotation angles at points (or along edges) where two
boundary faces meet, $\upsilon_{_1}$ is the combination of these
rotation angles with the rotation of the edges, and \(\tau\) is the
measure of the size of wedges at corners where three boundary faces of
\(\mM\) meet.
\end{theorem}

If the manifold $\mM$ is a single point (i.e., $d = 1$), the result
(\ref{eq:wts}) simplifies to $0.5 \, \mbP \left( \chi_1^2 \ge c
\right)$. If $\mM$ is one-dimensional (i.e., $d = 2$), the result
reduces to $(\kappa_{_0}/2 \, \pi) \, \mbP \left( \chi_2^2 \ge c
\right) + (\ell_0/4) \, \mbP \left( \chi_1^2 \ge c \right)$, where
$\kappa_{_0}$ is the length of $\mM$ and $\ell_0$ is the number of
boundary caps which equals 2. This last result is same as that
obtained by \cite{lin:97}; however, we provide an explicit formula for
$\kappa_{_0}$ based on the parametric representation of $\mM$ which is
derived in the next section.

\begin{remark}
For convex manifold $\mM$, $\kappa_{_2} = 0$ in the asymptotic
expansion (\ref{eq:wts}). In the case of order restricted
alternatives, the manifold $\mM$ is a high-dimensional tetrahedron
whose corners correspond to the constraints imposed on $\mbg$ under
$\mH_1$; hence $\upsilon_{_0} \neq 0$. Also, $\ell_2 = 0$ in the case
of order restricted testing problem. 
\end{remark}

Owing to Theorem \ref{th:sig3d}, the weights in Theorem 2.1 of
\cite{take:97} for the independent data case can be determined using
the volume-of-tube formula.  When the critical geodesic radius of the
manifold \citep{naiman:90} is greater than or equal to $\pi/2$, all of
the coefficients in (\ref{eq:wts}) are nonnegative leading to a finite
mixture of chi-squared distribution or $\overline{\chi}^2$
distribution. The critical radius is greater than equal to $\pi/2$ if
and only if the smallest cone containing the manifold is convex. Lemma
2.1 of \cite{take:02} provides a formula for computing the critical
radius. 

\vspace{0.2in}
\subsection{\em Explicit expressions for the geometric constants}
\vspace{-.1in}
We derive explicit expressions (in suitable forms for computation) for
the geometric constants in (\ref{eq:wts}) using the representation of
$\mT(\mbn)$, defined in (\ref{eq:Tn}), and its derivatives. Our main
goal is to reduce the evaluation of the constants to integrals over
appropriate parts of $\mM$.

The profound result of Gauss-Bonnet theorem \citep{do:76,milman:77}
connecting curvatures of manifolds with the Euler-Poincar\`e
characteristic \citep{worsley2:95,worsley1:95,adler:04} can be
employed to find some of the geometric constants appearing in Theorem
\ref{th:sig3d}. When $\mM$ is two-dimensional (i.e., $d = 3$), the
number of pieces contributing to $\mM$ minus the number of holes
equals $\mathcal{E}$. In particular, $\kappa_{_2} + \ell_1 +
\upsilon_{_0} = 2 \, \pi \mathcal{E} - \kappa_{_0}$ which eliminates
the need to compute $\kappa_{_2}, \ell_1$ and $\upsilon{_0}$
directly. \\
\begin{remark}
\cite{lin:97} assume that the cone is convex and smooth; hence no
corners (i.e., $\upsilon_{_0} = 0$). Both their Theorem 3.1 and
the result for $d = 3$ presented in Section 4 of \cite{lin:97} become
special cases of our general result established in Theorem
\ref{th:sig3d}.   
\end{remark} 

In the parametric representation of $\mM$, the function
$\mT(\mbn)$ has an embedded constraint $\|\bH \, \mbn\| = 1$ for a
$d$-dimensional vector $\mbn \in \mN$. This embedded constraint means
that the manifold $\mM$ is of dimension $(d - 1)$. Hence, we can
represent in terms of a $(d - 1)$-dimensional parameter vector
$\mbrh$ with $\mbn \equiv \mbn(\mbrh)$ and express $\mbT(\mbrh) =
\mT[\mbn(\mbrh)]$. For example, such a transformation can be carried
out using the polar co-ordinates. Denote the domain of $\mbrh$ by
$\mN^{\star}$. We express $\mbT(\mbrh) = [\mbT^1(\rho_{_1}, \ldots,
\rho_{_{d - 1}}), \ldots, \mbT^{\, r}(\rho_{_1}, \ldots, \rho_{_{d -
1}})]$ for the parametric representation of $\mM$. \\
\begin{assumption}
\label{as10}
The transformation $\mbT(\mbrh)$ is one-to-one and each $\mbT^l$ $(l =
1, \ldots, r)$ is twice continuously differentiable on $\mN^{\star} 
\subset \Re^{(d - 1)}$.   
\end{assumption}
\\

The following expressions are derived under \ref{as10}.  Define an $[r
\times (d - 1)]$ matrix $\mbQ(\mbrh) := [\mbT_1(\mbrh) \cdots \mbT_{d
- 1}(\mbrh)]$, where $\mbT_k(\mbrh) = \partial \, \mbT(\mbrh)/\partial
\, \rho_{_k}$ for $k = 1, \ldots, (d - 1)$ and $\mbT_1(\mbrh),
\ldots$, $\mbT_{d - 1}(\mbrh)$ are the column vectors.

The volume of the manifold $\mM$ is expressed as 
\bes
  \kappa_{_0} = \int_{\mbrh \, \in \, \mN^{\star}} \; \det \left[
  \mbQ^T(\mbrh) \; \mbQ(\mbrh) \right]^{1/2} \; d \, \mbrh. 
\ees
Finding $\ell_2$ is essentially similar to that of finding $\kappa_2$
by simply treating each face of the boundary as a new
manifold. Therefore, we describe the method to find $\kappa_{_2}$,
although it is zero for convex manifolds. In the case of testing under
order restricted alternatives, the constant $\ell_2 = 0$. The measure
of curvature of $\mM$ is expressed as
\bes
   \kappa_{_2} = \int_{\mbrh \, \in \, \mN^{\star}} \; \frac{1}{2}
   \left[ \Upsilon(\mbrh) - (d - 1)(d - 2) \right] \;
   \det \left[ \mbQ^T(\mbrh) \; \mbQ(\mbrh) \right]^{1/2} 
   \; d \, \mbrh, 
\ees
where $\Upsilon(\mbrh) = 2 \sum_{k = 2}^{d - 1} \, \sum_{l = 1}^{k -
1} \left( \mbv^T_{kk} \, \mbv_{ll} - \mbv_{kl}^T \, \mbv_{lk} \right)$
with  
\bes
  \mbv^T_{lk} = \b e_l^T \left[ \mbQ^T(\mbrh) \, \mbQ(\mbrh)
  \right]^{-1} \, \left[ \frac{\partial}{\partial \, \rho_{_k}}
  \mbQ(\mbrh) \right] \, \left[ \bI_{d - 1} - \mbQ(\mbrh) \, \left\{
  \mbQ^T(\mbrh) \, \mbQ(\mbrh) \right\}^{-1} \, \mbQ^T(\mbrh) \right] 
\ees
and $\b e_l$ as the basis vector for $\Re^{(d - 1)}$. 

The volume of $\partial \, \mM$, the boundary of $\mM$, is $\ell_0$
while $\ell_1$ measures the curvature of $\partial \, \mM$; both of
these need to be determined for each face of the boundary. For
instance, on the face where $\rho_{_{d - 1}} = 1$, define
$\mbQ_{\dag}(\mbrh) := [\mbT_1(\mbrh) \cdots \mbT_{d - 2}(\mbrh)]$ and
$\bB_{d - 1}(\mbrh) := \varphi \left[ \bI_{d - 1} - \mbQ(\mbrh) \,
\left\{\mbQ^T(\mbrh) \; \mbQ(\mbrh) \right\}^{-1} \; \mbQ^T(\mbrh)
\right] \, \mbT_{d - 1}(\mbrh)$, where $\varphi$ is a normalizing
constant. The geometric constants $\ell_0$ and $\ell_1$ can be
determined via 
\bes
\ell_0 = \int_{\mbrh \, \in \, \partial \, \mM} \, \det\left[
\mbQ_{\dag}^T(\mbrh) \, \mbQ_{\dag}(\mbrh) \right]^{1/2} \, d \, \mbrh
\ees
and 
\bes
   \ell_1 =  \int_{\mbrh \, \in \, \partial \, \mM} \ell_1(\mbrh) \;
   \det \left[ \mbQ_{\dag}^T(\mbrh) \; \mbQ_{\dag}(\mbrh)
   \right]^{1/2} \, d \, \mbrh,  
\ees
where 
\bes
  \ell_1(\mbrh) = - \sum_{k = 1}^{d - 2} \b e_k^T \left[
  \mbQ_{\dag}^T(\mbrh) \; \mbQ_{\dag}(\mbrh) \right]^{-1} \;
  \left[\frac{\partial } {\partial \, \rho_{_k}} \mbQ_{\dag}^T(\mbrh)
  \right]\, \bB_{d - 1}(\mbrh).
\ees
Similarly, 
\bes
   \upsilon_{_1} =  \int_{\mbrh \, \in \, \partial^2 \mM}
   \upsilon_{_1}(\mbrh) \; \det \left[ \mbQ_{\dag}^T(\mbrh) \;
   \mbQ_{\dag}(\mbrh) \right]^{1/2} \, d \, \mbrh,  
\ees
where $\partial^2 \, \mM$ is the region or corner at which two
boundary faces of $\mM$ meet,
\bes
  \upsilon_{_1}(\mbrh) = - \sum_{k = 1}^{d - 3} \b e_k^T \left[
  \mbQ_{\dag}^T(\mbrh) \; \mbQ_{\dag}(\mbrh) \right]^{-1} \;
  \left[\frac{\partial } {\partial \, \rho_{_k}} \mbQ_{\dag}^T(\mbrh)
  \right]\, \left[ \bB_{d - 2}(\mbrh) + \bB_{d - 1}(\mbrh)
  \right] \, \tan \left[ \frac{\phi(\mbrh)}{2} \right]
\ees
and $\phi(\mbrh)$ is the angle between $\bB_{d - 2}(\mbrh)$ and
$\bB_{d - 1}(\mbrh)$. Since $\mM$ is of dimension $(d - 1)$,
$\partial \, \mM$ and $\partial^2 \, \mM$ are of dimensions $(d - 2)$
and $(d - 3)$, respectively. Lastly, we consider the edges where two
boundary faces meet. If we consider the edge where $\rho_{_{d - 2}} =
\rho_{_{d - 1}} = 1$ and define $\mbQ_{\ddag}(\mbrh) := [\mbT_1(\mbrh)
\cdots \mbT_{d - 3}(\mbrh)]$, then
\bes
   \upsilon_{_0} = \int_{\mbrh \, \in \, \partial^2 \mM}
   \upsilon_{_0}(\mbrh) \; \det \left[ \mbQ_{\ddag}^T(\mbrh) \;
    \mbQ_{\ddag}(\mbrh) \right]^{1/2}  \, d \, \mbrh,
\ees
where $\upsilon_{_0}(\mbrh) = \mbox{arc} \, \cos[ \ip{ \bB_{d -
2}(\mbrh)}{ \bB_{d - 1}(\mbrh)}]$. 

In order to find $\tau$, we need to calculate the area of the
spherical triangle which is achieved by the Euler's formula: area of
the triangle equals $(\phi_1 + \phi_2 + \phi_3 - \pi)$, where $\phi_1,
\phi_2$ and $\phi_3$ are the three internal angles of the
triangle. \cite{loader2:04} describe the method of determining these
angles by first finding the vectors defining the corners of the
triangles. 

\section{Power Under a Sequence of Local Alternatives} 
\label{sec-lalt}
In this section, we derive an asymptotic lower bound for the power of
the GQS statistic under a sequence of local alternatives in
$\mK$. This plays an important role in comparing the result with a
test against the unrestricted alternative. To the best of the author's
knowledge, a lower bound has not been established in the literature
even for independent data; hence it would be an interesting one to
derive. 

From the parameterization defined in Section \ref{sec-prep}, we can
express $\mbg = \mbg_1 + \bP \, \mbn$ such that $\mbg_1 \in \mV$ and
$\mbn \in \mN$. Following the hypothesis (\ref{hypn}), we consider a
sequence of local alternatives of the form
\bea
  \label{local}
  \mbn_{_N} =  \frac{\mbn^{\star}}{\sqrt{N}} \quad \mbox{for a fixed
  vector} \quad \mbn^{\star} \in \mN \subset \Re^d \quad \mbox{and}
  \quad N = 1, 2, \ldots. 
\eea
From the derivation of $\mbr_2$ in (\ref{eq:theta2}), the relation
$\mbr_2 = \sqrt{N} \, \bH \, \mbn_{_N} = \bH \, \mbn^{\star}$ holds
under the sequence of local alternatives (\ref{local}).

As a first step, we define a statistic for testing $\mH_0\!: \mbg
\in \mV$ against the unrestricted alternative $\mH_2\!: \mbg \in
\mG$ as $\S_N^{\star} := \left[ \mQ_N(\overline{\mbg}) - \mQ_N(\wmbg)
\right]$, where $\overline{\mbg}$ and $\wmbg$ are defined in Section
\ref{sec-asymp}. Using the arguments similar to Theorem \ref{thm12}
and the result (\ref{eq:mbr}), one can establish that $\wmbr
\rightsquigarrow \bZ^{\ddag} \sim N_r( \bH \, \mbn^{\star}, \bI_r)$
under the sequence of local alternatives. 

The arguments given in \cite{robert:88} and \cite{pilla:05} yield the
following result. 
\begin{theorem}
\label{lemm8} 
The asymptotic local power of the unrestricted test statistic
$\S_N^{\star}$ for a sequence of alternatives (\ref{local}) is   
\bea
  \label{eq:apower}
  \underset{N \rightarrow \infty}{\lim} \mbP \left( \S_N^{\star}
  \geq b_1 \right) =  \mbP \left[ \chi^2_{r} (\delta^{2}) \geq b_1
  \right], 
\eea
where $b_1 > 0$ is a constant, $\delta = \|\bH \, \mbn^{\star}\|$ and
$\chi^2_{r} (\delta^{2})$ is the chi-square distribution with a
non-centrality parameter $\delta^{2}$ and with $r$ degrees of
freedom. 
\end{theorem}

The above result is equivalent to Theorem 7 of \cite{pilla:05};
however, here the non-centrality parameter $\delta$ is represented
in terms of $\bH$. Theorem \ref{lemm8} yields an exact local power and
the next one establishes a lower bound for $\S_N^{\star}$.
\begin{theorem}
A lower bound for the asymptotic power of $\S_N^{\star}$, under a
sequence of alternatives defined in (\ref{local}), is $\underset{N
\rightarrow \infty}{\lim} \mbP \left(\S_N^{\star} \geq b_1 \right) \geq
1 - \Phi \left( \sqrt{b_1} - \delta \right)$.
\end{theorem}

It is worth noting that finding the asymptotic power for $\S_N$ under
the sequence of local alternatives in $\mK$ is hard and it does not
have a simple weighted non-central chi-squared distribution with a
pre-specified non-centrality parameter. The following result gives an
asymptotic lower bound for the power of $\S_N$ under a sequence of
local alternatives (\ref{local}). It demonstrates that $\S_N$ under
restricted alternatives (i.e., testing for $\mH_0$ against $\mH_1$) is
locally more powerful than $\S_N^{\star}$ under no restriction (i.e.,
testing for $\mH_0$ against $\mH_2$).
\begin{theorem} 
\label{thm-lalt}
A lower bound for the asymptotic power of $\S_N$ for a sequence
of alternatives (\ref{local}) is $\underset{N \rightarrow
\infty}{\lim} \mbP \left( \S_N \geq b_2 \right) \geq \mbP \left(
N(\delta, 1) \geq  \sqrt{b_2} \right) = 1 - \Phi \left(\sqrt{b_2} -
\delta \right)$, where $b_2 > 0$ is a constant, $\delta := \|\bH \,
\mbn^{\star}\|$ and $\Phi$ is the standard Gaussian cumulative
distribution.   
\end{theorem}
{\sc Proof.} Let $\b d = \underset{\b \bb \, \in \,
\mK}{\arg \, \min} \; \| \bZ^{\ddag} - \b \bb\|^2$, where $\bZ^{\ddag}
\sim N_r( \bH \, \mbn^{\star}, \bI_r)$. Consequently, $\b d =
\mP_{\mK} \, \bZ^{\ddag}$. Further let, $\bbL := \{c \,  \bH \,
\mbn^{\star}\!: c > 0\}$ and $\b d^{\star} = \underset{\b \a \, \in
\, \bbL}{\arg \, \min} \; \| \bZ^{\ddag} - \b \a\|^2$, then $\b
d^{\star} = \mP_{\bbL} \, \bZ^{\ddag}$. Since $\bbL \subset \mK$, it 
follows that $\|\b d\| \geq \|\b d^{\star}\|$. Equivalently, $\|
\mP_{\mK} \, \bZ^{\ddag} \| \geq \|\mP_{\bbL} \, \bZ^{\ddag} \|$.  

Let $\b \a = \bH \, \mbn^{\star}/\| \bH \, \mbn^{\star} \|$ so that
$\|\b \a\| = 1$. It is clear that $\|\b d^{\star}\| = \ip{\b
\a}{\bZ^{\ddag}}_+$. The proof of $\ip{\b \a}{\bZ^{\ddag}} \sim
N(\delta, 1)$ is presented next. We have $\mE \left[ \ip{\b
\a}{\bZ^{\ddag}} \right] = \ip{ \bH \, \mbn^{\star}}{ \bH \,
\mbn^{\star}}/\| \bH \, \mbn^{\star}\| = \delta$ and variance
$\mbV(\|\b d^{\star}\|) = \| \bH \, \mbn^{\star}\|^{-2} \; \mbV[(\bH
\, \mbn^{\star})^T \, \bZ^{\ddag}] = 1$, since $\bZ^{\ddag}  \sim N_r(
\bH \, \mbn^{\star}, \bI_r)$. Therefore, the relation $\ip{\b
\a}{\bZ^{\ddag}} \sim N(\delta, 1)$ holds. Furthermore, $\|\mP{_\mK} \,
\bZ^{\ddag}\| \geq \ip{\b \a}{\bZ^{\ddag}}$, since the right-hand side
is the length of the projection of $\bZ^{\ddag}$ onto $ \bH \,
\mbn^{\star} \in \mK$ and $\ip{\b \a}{\bZ^{\ddag}}$ is distributed as
$N(\delta, 1)$. The theorem claim follows from
(\ref{eq:local}). \hfill \rule{2mm}{2mm}  
\\
\begin{example} (Comparison of the local power of $\S_N^{\star}$ and
$\S_N$ for the order restricted testing problem). Define $\mbg_{N}
:= (\mu_{_{1, \, N}}, \mu_{_{2, \, N}}, \ldots, \mu_{_{m, \, N}},
\mbt_1, \ldots, \mbt_m)^T$ as a sequence of local parameter vectors
for $N = 1, 2, \ldots$. The local alternatives take the form
$\mu_{_{k, \, N}} = \mu_{_{k - 1, \, N}} + \epsilon_{_{k -
1}}/\sqrt{N}$ for $k = 2, \ldots, m$, where $\epsilon_{_1}, \ldots,
\epsilon_{_{m - 1}}$ are fixed negative constants. As $N \rightarrow
\infty$, the sequence of local alternatives approach the null
hypothesis $\mH_0\!: \mbg \in \mV$.

Consider three treatment groups (i.e., $m = 3$) leading to six
possible orderings, with each order corresponding to an arc on the
unit circle. Union of these six arcs comprises the unit circle
$\mS^1$. Due to the balanced design assumption, each of these arcs is
of the same length; therefore, the angle of the cone $\mK$ is $\phi =
\pi/3 = 60^{\circ}$. At the level of significance $\alpha = 0.05$, the
critical values corresponding to the two tests $\S_N$ and
$\S_N^{\star}$ are $b_2 = 3.820$ and $b_1 = 5.991$, respectively. Table
\ref{power} presents the asymptotic lower bounds on the local power
for the two tests. The table also presents the exact asymptotic local
power for $\S_N^{\star}$ obtained using the asymptotic formula
(\ref{eq:apower}). It is clear that except for $\delta = 0$, the
asymptotic local power of $\S_N$ is better than that of $\S_N^{\star}$,
in terms of both the lower bound and the exact power.
\end{example}
\vspace{-0.15in}
\begin{table}[htp]
\caption{Comparison of the local power of the tests for a given
$\delta$, the non-centrality parameter, when $m = 3$ and $\phi = \pi/3
= 60^{\circ}$.}
\label{power}
\centerline{\begin{tabular}{|l|cccccc|} \hline
Test & \multicolumn{6}{c|}{$\delta$} \\
  & 0 & 1 & 2 & 3 & 4 & 5 \\ \hline 
$\S_N$ Lower bound & 0.025 & 0.170 & 0.518 & 0.852 & 0.980 & 0.999 \\
$\S_N^{\star}$ Exact & 0.050 & 0.133 & 0.416 & 0.771 & 0.957 & 0.996
\\
$\S_N^{\star}$ Lower bound & 0.007 & 0.074 & 0.327 & 0.710 & 0.940 &
0.995 \\ \hline
\end{tabular}}
\end{table}
\section{Discussion}
\label{sec-discuss}
In this research, inferential theory is developed for the problem of
testing under convex cone alternatives for correlated data. Such a
problem occurs when interest lies in detecting ordering of treatment
effects, while simultaneously modeling relationships with other
covariates. The testing problem (\ref{eq:hyp0}) is also applicable to
the analysis of clustered multi-categorical data. In this framework, 
$\bY_{\!it} = (Y_{i1t}, \ldots, Y_{iKt})^{\, T}$ denotes the
$K$-categorical response on the $i$th observation in the cluster $t$,
where $Y_{ijt} = 1$ if category $j$ ($j = 1, \ldots, K$) is observed
and 0 otherwise.

We established that the GQS statistic is asymptotically equivalent to
the squared length of the projection of the standard Gaussian vector
onto an arbitrary convex cone with a nonempty interior. We further
derived the asymptotic null distribution of the GQS statistic under
convex cone alternatives for correlated data as a weighted chi-squared
distribution. The weights in the asymptotic distribution are the mixed
volumes of the convex cone and its polar cone which do not have
explicit expressions except in special cases. For non-polyhedral
cones, closed-form expressions for the weights are very complicated
and therefore; often a simulation approach is employed for computing
them [Section 3.5, \cite{silva:04}]. In this article, explicit
formulas are derived for the calculation of these weights using the
Hotelling-Weyl-Naiman volume-of-tube formula.

Furthermore, an asymptotic lower bound is derived for the power of the
test under a sequence of local alternatives in $\mK$ for correlated
data which establishes that the test under restricted alternatives is
more powerful than the test under no restriction. Note that
\cite{barlow:72} and \cite{robert:88} derive the asymptotic power only
under specified alternative hypothesis.  The current theory is
applicable to many practical problems of interest including testing
for a monotone regression function and for the analysis of clustered
multi-categorical data. 

\indent {\bf Acknowledgments.} The author is grateful to Catherine
Loader for many stimulating discussions and to the Associate Editor
for constructive comments. 

\bibliographystyle{apalike}
\bibliography{gcone}
\Line{\AOSaddress{
Department of Statistics \\ 
Case Western Reserve University\\
Cleveland, OH 44106 \\
Email: pilla@case.edu}\hfill
\AOSaddress{}}

\end{document}